\documentstyle[12pt]{article}
\def\hybrid{\topmargin 0pt      \oddsidemargin 0pt
        \headheight 0pt \headsep 0pt
        \textwidth 160true mm       
        \textheight 231true mm         
        \marginparwidth 0.0in
        \parskip 0pt plus 1pt   \jot = 1.5ex}
\catcode`\@=11
\def\marginnote#1{}
\newcount\hour
\newcount\minute
\newtoks\amorpm
\hour=\time\divide\hour by60 \minute=\time{\multiply\hour by60
\global\advance\minute by-\hour}
\edef\standardtime{{\ifnum\hour<12 \global\amorpm={am}%
        \else\global\amorpm={pm}\advance\hour by-12 \fi
        \ifnum\hour=0 \hour=12 \fi
        \number\hour:\ifnum\minute<10
        0\fi\number\minute\the\amorpm}}
\edef\militarytime{\number\hour:\ifnum\minute<10
0\fi\number\minute}

\def\draftlabel#1{{\@bsphack\if@filesw {\let\thepage\relax
   \xdef\@gtempa{\write\@auxout{\string
      \newlabel{#1}{{\@currentlabel}{\thepage}}}}}\@gtempa
   \if@nobreak \ifvmode\nobreak\fi\fi\fi\@esphack}
        \gdef\@eqnlabel{#1}}
\def\@eqnlabel{}
\def\@vacuum{}
\def\draftmarginnote#1{\marginpar{\raggedright\scriptsize\tt#1}}
\def\draft{\oddsidemargin -.5truein
        \def\@oddfoot{\sl preliminary draft \hfil
        \rm\thepage\hfil\sl\today\quad\militarytime}
        \let\@evenfoot\@oddfoot \overfullrule 3pt
        \let\label=\draftlabel
        \let\marginnote=\draftmarginnote
   \def\@eqnnum{(\theequation)\rlap{\kern\marginparsep\tt
   \@eqnlabel}%
\global\let\@eqnlabel\@vacuum}  }
\relax

\hybrid

\def\cp{{\bf CP}^n} \def\gggg{{\bf g}}
\title{Cotangent and tangent modules on quantum orbits}
\author{P.~Akueson,\\
D.~Gurevich,\\
ISTV, Universit\'e de Valenciennes\\
59304 Valenciennes, France}

\begin{document}

\maketitle

\begin{abstract}
Let $k(S^2_q)$ be the "coordinate ring" of a quantum sphere. We
introduce the cotangent module on the quantum sphere as a
one-sided $k(S^2_q)$-module and show that there is no Yang-Baxter
type operator converting it into a $k(S^2_q)$-bimodule which
would be a flatly deformed object w.r.t. its classical
counterpart. This implies non-flatness of any covariant
differential calculus on the quantum sphere making use of the
Leibniz rule. Also, we introduce the cotangent and tangent
modules on generic quantum orbits and discuss some related
problems of "braided geometry".
\end{abstract}



\def\to{\longrightarrow}
\def\bbbr{{\rm I\!R}}     
\def\bbbc{{\rm I\!\!\! C}}  
\def\bbbp{{\rm I\!P}}
\def\bbbk{{\rm I\!K}}
\def\uqs{ U_q(sl(2))}
\def\Usq{ U_q(sl(n))}
\def\uqsm{ U_q(sl(2))-Mod}
\def\ot2{\otimes 2}
\def\ovalp{\overline{\alpha}}
\def\ovb{\overline{\beta}}
\def\ovs{\overline{S}}
\def\ovc{\overline{C_q}}
\newtheorem{thm}{Th\'eor\`eme}[section]
\newtheorem{lemme}{Lemme}[section]
\newtheorem{corol}{Corollaire}[section]
\newtheorem{defin}{D\'efinition}[section]
\newtheorem{rem}{Remarque}[section]
\newtheorem{exem}{Exemple}[section]
\newtheorem{prop}{Proposition}[section]
\newcommand{\uq}{U_q(\gggg)}
\newcommand\lbc{[}
\newcommand\rbc{]}
\newcommand{\Pm}{P_-(t)}
\newcommand{\rr}{\rho^{\ot 2}}

\renewcommand{\theequation}{{\thesection}.{\arabic{equation}}}
\def\r#1{\mbox{(}\ref{#1}\mbox{)}}
\setcounter{equation}{0}
\def\ot{\otimes}
\def\kMq{k(M_q)}
\def\Uq{U_q({\bf g})}
\def\kSq{k(S^2_q)}
\def\oS{\overline S}
\def\kS{k(S^2)}
\def\uq{U_q(sl(2))}
\def\Uqs{U_q(sl(n))}
\def\uqsl{U_q(sl(2))}
\def\TSl{T^*_l(S^2_q)}
\def\TSr{T^*_r(S^2_q)}
\def\TSL{T^*_l(S^2)}
\def\TSR{T^*_r(S^2)}
\def\oTS{\overline{T^*(S^2)}}
\def\oTSq{\overline{T^*(S^2_q)}}
\def\TS{T^*(S^2)}
\def\TSq{T^*(S^2_q)}
\def\aam{A^{(m)}}
\def\aaaa{\cal A}
\def\wm{\wedge_-}
\def\wp{\wedge_+}
\def\wpm{\wedge_{\pm}}
\def\vn{V^{\ot n}}
\def\iin{I^{(n)}}
\def\inn{I^{n}}
\def\de{\delta}
\def\De{\Delta}
\def\RR{\bf R}
\def\CC{\bf C}
\def\cp{{\bf CP}^n}
\def\Mat{\rm Mat}
\def\Sym{\rm Sym}
\def\ahqc{A^c_{\h,\,q}}
\def\Om{\Omega}
\def\aqc{{\cal A}^c_{ 0\, q}}
\def\ahqc{{\cal A}^c_{\h\, q}}
\def\ac{{\cal A}^c_{0\, 1}}
\def\h{{\hbar}}
\def\ah{{\cal A}_{\hbar}}
\def\vv{V^{\ot 2}}
\def\thq{T(H_q)}
\def\am{A_{\mu}}
\def\ami{A_{\mu}^i}
\def\gq{{\gggg}_q}
\def\ggq{{\gggg}_q^{\ot 2}}
\def\sss{\cal S}
\def\Ob{Ob\,}
\def\h{{\hbar}}
\def\a{\alpha}
\def\na{\nabla}
\def\ve{\varepsilon}
\def\ep{\epsilon}
\def\bea{\begin{eqnarray}}
\def\eea{\end{eqnarray}}
\def\beq{\begin{equation}}          \def\bn{\beq}
\def\eeq{\end{equation}}            \def\ed{\eeq}
\def\nn{\nonumber}
\def\oH{{\overline H}}
\def\aaaa{\cal A}
\def\aaa{\cal A}
\def\End{{\rm End\, }}
\def\uEnd{\underline{\rm End}\,}
\def\Ker{{\rm Ker\, }}
\def\root{{\rm root\, }}
\def\diag{{\rm diag\, }}
\def\Im{{\rm Im\, }}
\def\Id{{\rm Id }}
\def\ra{{\rm rank\, }}
\def\Id{{\rm id\, }}
\def\Vect{{\rm Vect\,}}
\def\Hom{{\rm Hom\,}}
\def\uHom{\underline{\rm Hom}\,}
\def\tr{{\rm tr}}
\def\span{{\rm span}}
\def\det{{\rm det}}
\def\codet{{\rm codet}}
\def\dim{{\rm dim}\,}
\def\id{{\rm id}\,}
\def\span{{\rm span}\,}
\def\udim{\underline{\rm dim}\,}
\def\dem{{\rm det}^{-1}}
\def\Fun{{\rm Fun\,}}
\def\de{\delta}
\def\al{\alpha}
\def\cO{\cal O}
\renewcommand{\theequation}{{\thesection}.{\arabic{equation}}}
\setcounter{equation}{0}
\newtheorem{proposition}{Proposition}
\newtheorem{conjecture}{Conjecture}
\newtheorem{corollary}{Corollary}
\newtheorem{theorem}{Theorem}
\newtheorem{definition}{Definition}
\newtheorem{remark}{Remark}
\def\R{{\cal R}}
\def\n{{\bf n}}   \def\va{\varepsilon}
\def\so{s_{\omega_1+\omega_{n-1}}}

\section{Introduction}

Since the creation of super-theory it became clear that numerous
aspects of commutative algebra and usual geometry could be
generalized to the super-case. In particular, for any two
$Z_2$-graded one-sided (say, left) $A$-modules $M_1$ and $M_2$
over a super-commutative algebra $A$ their tensor product
$M_1\ot_A M_2$ is well-defined.

In the latter 80's it was recognized that many properties of
(super-)commutative algebra and geometry could be further
generalized onto objects related  to a Yang-Baxter (YB) operator,
i.e., a solution of the quantum YB equation $$ (S\ot \id)(\id \ot
S)(S\ot \id)=(\id \ot S)(S\ot \id)(\id \ot S), $$ $S$ being an
operator acting on $V^{\ot 2}$ where $V$ is a vector space. If
$S$ is an involutary  YB operator $(S^2=\id)$, the notion of an
$S$-commutative algebra can be introduced in a natural way. If
$A$ is such an algebra, the product $M_1\ot_A M_2$ of two
one-sided $A$-modules can be introduced by means of the operator
$S$ (under some natural conditions on it). Also, the operator $S$
plays the crucial role in a twisted or quantum (i.e. related to
an operator $S$) version of differential calculus. It can be used
for ordering "quantum functions" and differentials mixed in
virtue of the Leibniz rule.

Unfortunately, a straightforward application of this method to
algebras related to a non-involutary YB operator $S$ leads to
non-flat deformations\footnote{Let us recall that a deformation
$V_{\h}$ of a vector space $V$  where $\h$ is a formal parameter
is called flat if $V_{\h}/\h V_{\h}=V$ and $V_{\h}$ is isomorphic
to $V\ot k[[\h]]$  as $k[[\h]]$-module (the tensor product is
completed  in the $\h$-adic topology).}.

As was shown in \cite{Ar}, \cite{AAM} the  differential calculus
on the quantum algebras $$k_q(G), G=SL(n), SO(n), Sp(n)$$
 initiated in \cite{W1} and \cite{W2} and making use of the Leibniz rule
gives rise to non-flat deformations of the classical differential
algebras. However, this differential calculus plays the central
role in all known attempts to introduce a quantum  version of
gauge theory related to Drinfeld-Jimbo quantum groups (cf., i.e.
 \cite{BM}, \cite{HM} and the references therein). From
our viewpoint, non-flat deformations are somewhat pointless since
in this case classical objects are not limits of their quantum
counterparts.

Quantum orbits pertain to another type of algebras related to
quantum groups. The simplest example of such an orbit is the
quantum sphere $\kSq$ introduced in \cite{P1}. A version of
differential calculus on it  was suggested in \cite{P2}. However,
as follows from our results, on the quantum sphere no flatly
deformed $U_q(su(2))$-covariant differential calculus exists
which makes use of the Leibniz rule. An explanation  of this
phenomenon, given in the present paper, consists in the
following. It is possible to realize a flat deformation of the
cotangent module\footnote{We use this term for the module
corresponding to the tangent vector bundle in the framework of
the Serre-Swan approach, cf. \cite{Se}, \cite{Sw}.
 In a similar way we
consider other classical and quantum modules. (All modules are
assumed to be finitely generated.)} considered as a one-sided
$k(S^2)$-module.  However, the flatness of deformation breaks
down  when one tries to deform the tangent vector bundle
considered as a two-sided $k(S^2)$-module.

Thus, if we want deformed modules to be flatly deformed objects,
the use of one-sided modules on quantum varieties \footnote{By
abusing the language we call "quantum varieties" the
corresponding "coordinate rings".} are only relevant in general.

Note, that some one-sided modules which are q-analogues of line
bundles on "quantum generic orbits" were constructed in \cite{GS}
in the spirit of the Serre-Swan approach. Here we introduce the
cotangent and tangent modules on the same orbits. Also,  we
discuss the problem of an operator meaning of the tangent
modules, i.e. that of representing them by "braided vector
fields" (note that in the case of the quantum sphere this problem
was solved in \cite{A}).

In the sequel we prefer to use the quantum group (QG) $\Uqs$
instead of this $U_q(su(n))$ since we disregard any involution in
quantum algebras in question. The basic field $k$ is assumed to
be $\bbbr $ or $\bbbc$. Throughout the whole of the paper the
parameter $q$ is assumed to be generic.

{\bf Acknowledgment} One of the authors (D.G.) is supported by
the grant PICS-608.

\section{Cotangent $\kSq$-module}

First, let us introduce the quantum sphere in a form appropriate
for our goals. Let $V$ be a three dimensional vector space being
a left $\uq$-module.  The action of the quantum group $\uq$ can
be extended to any tensor power of $V$ via the coproduct. Let us
decompose the space $\vv$ into a direct sum of irreducible
$\uq$-modules:

\beq \vv=V_0 \oplus \, V_1 \oplus \, V_2 \label{dec} \eeq where
$V_i$ stands for the spin $i$ $\uq$-module. Let $v_0$ be a
generator of the one-dimensional component $V_0$.

Let us introduce quantum sphere by \beq \kSq=T(V)/\{V_1,\,
v_0-c\},\,\, c\in k,\,\, c\not=0 \label{qs} \eeq where $T(V)$
stands for the free tensor algebra of the space $V$ and $\{I\}$
stands for the ideal generated by a family $I\subset T(V)$. This
algebra is a particular (in some sense q-commutative) case of the
Podles' sphere introduced in \cite{P1}. (We disregard involution
operators in this algebra. So, in fact we do not distinguish
quantum sphere and quantum hyperboloid.)

Remark, that in definition \r{qs} of the quantum sphere we do not
use any coordinate form of the quantum sphere. In a similar
non-coordinate way we introduce the cotangent module on it.

Let $V'$ be the space isomorphic to $V$ as $\uq$-modules but
spanned by the differentials $d\,x,\,\,x\in V$. Let $V'\ot\kSq$
be the free finitely generated right $\kSq$-module and  its
submodule $M_r$ generated by  $(V'\ot V)_0$. Thus, we have $$
M_r=\Im\,\mu^{23}((V'\ot V)_0\ot \kSq). $$ Hereafter $(V'\ot
V)_i$ stands for the spin $i$ component in the product $V'\ot
V$,  $\mu$ stands for the product in the algebra $\kSq$, and
$\mu^{ij}$ is the operator $\mu$ applied to the i-th and j-th
factors.  We call {\em a right cotangent module} on the quantum
sphere the following quotient $$ \TSr=(V'\ot\kSq)/M_r. $$ We
define {\em a left cotangent module} $\TSl$ on the quantum sphere
in the same way as the quotient of the left $\kSq$-module
$\kSq\ot V'$ over its submodule $$ M_l= \mu^{12}(\kSq\ot (V\ot
V')_0). $$ (Hereafter we omit the symbol $\,\Im$.)

These modules are quantum analogues (respectively, right and left
ones) of the cotangent vector bundle over quantum sphere (or
quantum hyperboloid) realized in the spirit of the Serre-Swan
approach.

As was shown in \cite{AG}, the left cotangent module $\TSl$ is a
flat deformation of its classical counterpart (in fact, the proof
consists in showing that the classical and quantum objects are
built from the same, respectively, $U(sl(2))$- and
$\uq$-irreducible components).  By the same reason the right
cotangent module $\TSr$ is a flatly deformed object.

Now,  define a two-sided cotangent module (in the sequel called
{\em cotangent bimodule}) on the quantum sphere. Let us set $$
\oTSq=(\kSq\ot V'\ot \kSq)/(M_l\ot\kSq+\kSq\ot M_r). $$ This
$\kSq$-bimodule is much bigger than  one-sided one even in the
classical case $(q=1)$ because so far we do not have any rule for
transposing the "quantum functions" and differentials (i.e.,
elements of $V'$). In what follows we omit the subscription $q$
if $q=1$.

First, let us consider the classical case in details. In order to
reduce this bimodule to the seize of the one-sided one we should
define a commutation rule between elements of the algebra and
those of $V'$. In the classical case it is always done by the
flip. Namely, we set \beq \TS=\oTS/\{a\ot\,v-v\ot\,a\} \qquad
a\in k(S^2),\, v\in V'. \label{mod} \eeq

It is not difficult to see that there exists a one-to-one
correspondence between the one-sided (say, right) $\kS$-module
$\TSR$ and the bimodule $\TS$. Indeed, modulo the denominator of
\r{mod} any element of the $k(S^2)$-bimodule $k(S^2)\ot V'\ot
k(S^2)$ can be reduced to an element of the right module $V'\ot
k(S^2)$. Thus, we have a map $$ \rho: \,k(S^2)\ot V'\ot k(S^2)\to
V'\ot k(S^2). $$ The following inclusion is clear $$
\rho\,(M_l\ot k(S^2)+k(S^2)\ot M_r) \subset M_r. $$ This implies
that the map $\rho$ sends the two-sided module $\TS$ into $\TSR$.
Moreover, it is isomorphism of linear spaces.

An analogous construction for algebras related to an involutary
YB operator $S$ can be introduced in a similar way. In this case
the denominator of formula \r{mod} should be replaced by $\{a\ot
v-S(a\ot v)\}$. However, if we want to realize a similar approach
for the algebra $\kSq$ (or for other algebras related to
non-involutary YB operators $S$) it is not clear what should be a
proper analogue of the denominator in \r{mod}.

Let $\oS: V\ot V'\to V'\ot V$ be any $\uq$-covariant invertible
operator (called in the sequel {\em a transposition}). Let us
replace the denominator in \r{mod} by  $\{a\ot v-\oS (a\ot v)\}$
(we also assume that $S^2$ in the numerator is replaced by
$S^2_q$). The problem consists in finding all transpositions
$\oS$ such that the corresponding quotient denoted $\TSq$ would
be a flat deformation of the quotient \r{mod}.

It is evident that in order to give rise to a flatly deformed
object a transposition $\oS$ should preserve the ideal of the
formula \r{qs} and  take the submodule $M_l$ into $M_r$.
Otherwise, by passing to the right $\kSq$-module $\TSR$ we would
get some supplementary relations in it what would lead to a
collapse  of the final object. Essentially, this means that the
map $\rho$ takes the $\kSq$-bimodule $\oTSq$ onto some proper
quotient of the right $\kSq$-module $\TSR$.

As we will see in the next section the only transpositions
preserving the ideal in \r{qs} are $\pm S^{\pm 1}$ where $S$ is
the YB operator coming from $U_q(sl(2))$.  Then we will show that
even these operators do not take the submodule $M_l$ into $M_r$.
So, any transposition $\oS$ leads to the collapse mentioned above.

\section{Non-existence of a flat deformation}

Now we represent the quantum sphere in a more explicit
(coordinate) form. Let us fix the base $(u,v,w)$ in the space $V$
with the following action of the QG $\uq$ $$
\begin{array}{cclccclcccl}
X.u&=&0,& \  & X.v&=&-(q+q^{-1})u,& \  & X.w&=&v,\\
Y.u&=&-v,& &Y.v&=&(q+q^{-1})w,& & Y.w&=&0,\\
H.u& = &2u,& & H.v&=&0, & & H.w&=&-2w.
\end{array}
$$ Hereafter $X, H, Y$ are the standard generators of the QG
$\uqs$ (cf. \cite{CP}).

Note, that the QG $\uqs$ acts on the space $V'$ in the same way
(we should only replace the generators $(u,v,w)$ by  $(du,dv,dw)$
in the formulae above).

Thus, the spaces $V_i, \, i=0,1,2$ being irreducible
$\uq$-modules are as follows $$
V_0=\span(v_0),\,\,v_0=(q^3+q)uw+v^2+(q+q^{-1})wu, $$ $$
V_1=\span(q^2uv-vu,\,(q^3+q)(uw-wu)+(1-q^2)v^2,\, wv-q^2vw), $$ $$
V_2=\span(u^2,\,uv+q^2vu,\, uw-qv^2+q^4wu,\, vw+q^2wv,\, w^2) $$
(the sign $\ot$ is systematically omitted).

It is well known that the YB operator $S$ being restricted onto
each component becomes scalar. Namely, $$
S\vert_{V_0}=q^{-4}\,\id,\,\,S\vert_{V_1}=-q^{-2}\,\id,\,\,
S\vert_{V_2}=q^{2}\,\id $$ (cf. f.e. \cite{G}). This implies that
being applied to the product $V\ot V'$ the operator $S$ acts as
follows
\begin{equation}
S(udu)=\alpha\,duu \label{eq1}
\end{equation}
\begin{equation}
S(q^2 udv -vdu)=\beta\,(q^2duv-dvu) \label{eq2}
\end{equation}
\begin{equation}
S((q^3+q)udw + vdv + (q + q^{-1})wdu)=\gamma ((q^3+q)duw + dvv +
(q + q^{-1})dwu) \label{eq3}
\end{equation}
with $\alpha= q^2,\, \beta=-q^{-2},\,\gamma=q^{-4}$ and similarly
for other elements of each component.

Now  consider an arbitrary invertible transposition $\oS:V\ot
V'\to V'\ot V$ being a $\uq$-morphism (we do not require it to be
a YB operator). It is given by the same formulae \r{eq1}-\r{eq3}
(and all their descendants) but with arbitrary non-trivial
$\alpha, \beta$ and $\gamma$. By applying such an operator  many
times we can transform any element from $V^{\ot k}\ot V'$ into
that from $V'\ot V^{\ot k}$. Let us note that the operator
$\oS=S$ has the following remarkable property (this property is
also valid for the operator $S^{-1}$).

\begin{proposition}
We have $$ S\mu^{12}=\mu^{23}S^{12}S^{23} \quad {\rm and} \quad
S\mu^{23}=\mu^{12}S^{23}S^{12} $$.
\end{proposition}

{\bf Proof} The defining relations of the algebra $\kSq$  are
coordinated with the action of the QG $\uqs$ in the following
sense $$ X.\mu(a\otimes b) = \mu \Delta(X).(a\otimes b)\quad
\forall X\in \uqs ,\, a,\, b \in \kSq. $$ Let ${\cal R}$ be the
universal quantum R-matrix corresponding to the QG $\uqs$. It
satisfies the relations $$ \Delta^{12}{\cal R}={\cal R}^{13}{\cal
R}^{23}\,\,\mbox{and}\,\, \Delta^{23}{\cal R}={\cal R}^{13}{\cal
R}^{12} $$ (this means that the QG in question is
quasitriangular).

Thus, for any $a,\, b, \, c\in \kSq$ we have (hereafter $\sigma$
is the usual flip)
\begin{eqnarray*}
S\mu^{12}(a\otimes b\otimes c)&=&\sigma {\cal R}\mu^{12}(a\otimes
b\otimes c)\\
&=&\sigma ({\cal R}_1(\mu^{12}(a\otimes b))\otimes {\cal R}_2
c)\\
&=&\sigma (\mu^{12}\Delta({\cal R}_1)(a\otimes b)\ot {\cal R}_2
c) = \sigma \mu^{12}(\Delta^{12}{\cal R}) (a\otimes b\otimes c).
\end{eqnarray*}
Here we use the Sweedler's notation and apply the components of
$\cal R$ to the elements $a,\,b,\, c$ w.r.t. the action $\uqs$ on
the algebra $\kSq$. Moreover, we use the relation $S=\sigma {\cal
R}$.

The following chain of identities completes the proof of the
first relation of proposition (the second one can be verified in
a similar way)
\begin{eqnarray*}
\sigma \mu^{12}\Delta^{12}{\cal
R}&=&\mu^{23}\sigma^{12}\sigma^{23}\Delta^{12}{\cal R}\\
&=&\mu^{23}\sigma^{12}\sigma^{23}{\cal R}^{13}{\cal
R}^{23}=\mu^{23}\sigma^{12}\sigma^{23} {\cal
R}^{13}\sigma^{23}\sigma^{23}
{\cal R}^{23}\\
&=&\mu^{23}\sigma^{12}\sigma^{23}{\cal
R}^{13}\sigma^{23}S^{23}=\mu^{23}\sigma^{12}
{\cal R}^{12}S^{23}\\
&=&\mu^{23}S^{12}S^{23}.
\end{eqnarray*}

It is evident that the above proposition is still valid if we
replace the operator $S$ by $S^{-1}$ in the formulae above.
However, if we replace the operator $S$ by any other
transposition $\oS$ the identities from this proposition become
broken in virtue of the following.

\begin{proposition}
The only operators $\oS$ such that $$ \oS(v_0\,v')=v'\,
v_0,\,\,\forall v'\in V'\,\,{\rm and} \,\, \oS (V_1 \otimes
V')\subset V'\otimes V_1 $$ are $\pm S$  and $\pm S^{-1}$.
\end{proposition}

{\bf Proof} In the sequel we represent an arbitrary transposition
as $$ \oS=xP_{0}+yP_{1}+zP_{2},\,\, x,\, y,\, z\in k $$ where the
operators
\begin{equation}
P_{i}:V\otimes V^{\prime }\rightarrow (V^{\prime }\otimes
V)_{i},\,\, i\in \left\{ 0,1,2\right\}
\end{equation}
become the projectors  $V^{\otimes 2}\rightarrow V_{i}$ if we
identify $V$ and $V'$.  For $\oS=S$  we have  $x=\gamma =q^{-4},$
\ $y=\beta =-q^{-2},$ \ $z=\alpha =q^{2}.$  In what follows we
need the images of some elements under the action of the
operators $P_i$:
\begin{eqnarray*}
P_{0}(vdv) &=&\alpha _{1}d\overline{v}_{0},\,\,P_{0}(udw)=\alpha
_{1}^{\prime }d\overline{v}_{0},\,\,P_{2}(udu)=duu, \\
P_{2}(udv) &=&\beta ^{\prime
}(duv+q^{2}dvu),P_{2}(vdu)=\overline{\beta }%
^{\prime }(duv+q^{2}dvu), \\
P_{2}(vdv) &=&\gamma _{1}(duw-qdvv+q^{4}dwu), \\
P_{1}(vdv) &=&\beta _{1}[(q^{3}+q)(duw-dwu)+(1-q^{2})dvv], \\
P_{1}(vdu) &=&\overline{\alpha }^{\prime }(q^{2}duv-dvu),\,\,
P_{1}(udv)=\alpha ^{\prime }(q^{2}duv-dvu), \\
dv_{0} &=&(q^{3}+q)udw+vdv+(q+q^{-1})wdu,\,\,d\overline{v}%
_{0}=(q^{3}+q)duw+ \\
&&dvv+(q+q^{-1})dwu,\\
P_{1}^{12}(ud\overline{v}_{0}) &=&\alpha ^{\prime
}(q^{2}duv-dvu)v+\beta
_{1}^{\prime }(q+q^{-1})[(q^{3}+q)(duw-dwu)+ \\
&&(1-q^{2})dvv]u, \\
P_{2}^{12}(ud\overline{v}_{0}) &=&(q^{3}+q)duuw+\beta ^{\prime
}(duv+q^{2}dvu)v+\gamma _{1}^{\prime }(q+q^{-1})(duw- \\
&&qdvv+q^{4}dwu)u, \\
P_{0}^{12}(ud\overline{v}_{0}) &=&(q+q^{-1})\alpha _{1}^{\prime
}d\overline{v%
}_{0}u,
\end{eqnarray*}
where
\begin{eqnarray*}
\beta ^{\prime } &=&(1+q^{4})^{-1},\,\,\alpha ^{\prime
}=q^{2}\beta ^{\prime },\,\,\overline{\alpha }^{\prime }=-\beta
^{\prime },\,\, \overline{\beta }^{\prime }=q^{2}\beta ^{\prime
},\gamma _{1}^{\prime
\prime }=q^{4}\gamma _{1}^{\prime }, \\
\beta _{2} &=&\beta ^{\prime },\,\,\alpha _{2}=-\alpha ^{\prime
},%
\,\,\alpha _{2}^{\prime }=\beta ^{\prime },\,\,\beta _{2}^{\prime
}=q^{2}\alpha _{2}^{\prime },\,\,\alpha _{1}^{\prime \prime
}=q^{-2}\alpha _{1}^{\prime }, \\
\alpha _{1} &=&q^{2}(1+q^{2}+q^{4})^{-1},\,\,\beta
_{1}=(1-q^{2})\beta ^{\prime },\,\,\beta _{1}^{\prime \prime
}=-\,\,\beta
_{1}^{\prime }, \\
\gamma _{1} &=&-q(1+q^{2})^{2}(1+q^{2}+q^{4})^{-1}\beta ^{\prime
},\,\,
\gamma _{1}^{\prime }=(1+q^{2}+q^{4})^{-1}\beta ^{\prime }, \\
\alpha _{1}^{\prime }
&=&q^{3}(1+q^{2})^{-1}(1+q^{2}+q^{4})^{-1},\,\, \beta
_{1}^{\prime }=q(1+q^{2})^{-1}\beta ^{\prime }.
\end{eqnarray*}

On applying the transposition $\oS$ to the element
$(q^{2}uv-vu)du$ we get
\begin{eqnarray*}
\overline{S}((q^{2}uv-vu)du) &=&\overline{S}^{12}\overline{S}%
^{23}((q^{2}uv-vu)du) \\
&=&q^{2}\overline{S}^{12}\overline{S}^{23}(uvdu)-\overline{S}^{12
}\overline{S%
}^{23}(vudu).
\end{eqnarray*}

By a straightforward but tedious computations with the use of the
formulae above we get the following result for the coefficient at
the element $dv\, uu$ in the image above
\begin{equation}
-q^{4}\beta ^{\prime 2}(\alpha ^{2}+(q^{4}+q^{-4})\alpha \beta
+\beta^{2}).
\end{equation}

The condition \ $\overline{S}(V_{1}\otimes V^{\prime })\subset
(V^{\prime }\otimes V_{1})$ implies
\begin{equation}
\alpha ^{2}+(q^{4}+q^{-4})\alpha \beta +\beta ^{2}=0.
\end{equation}

This equation has two solutions (in the projective sense)
\begin{equation}
\alpha =-q^{4}\beta\quad {\rm and}\quad \alpha =-q^{-4}\beta .
\label{tri}
\end{equation}

Let us remark that the first (resp., second) solution is
satisfied by the operator \ $cS$ (resp. $cS^{-1}$) \ with an
arbitrary $c\neq 0.$ Then, the operator  $\overline{S}$ in
general can be represented as follows
\begin{equation}
\overline{S}=cS^{\pm 1}+\delta P_{0},\,\, \delta\in k.
\end{equation}

Now it remains to show that  $\delta =0$ and  $c=\pm 1$. Let us
do it for $S$ (the \ $S^{-1}$ case is analogous). From the above
form of $\oS$  we have
\begin{equation}
\overline{S}^{12}\overline{S}^{23}=c^{2}S^{12}S^{23}+c\delta
(S^{12}P_{0}^{23}+P_{0}^{12}S^{23})+\delta
^{2}P_{0}^{12}P_{0}^{23}.
\end{equation}

Let us consider the image  of the element  $(q^{2}uv-vu)dv$
w.r.t. the action of  the transposition  $\overline{S}$. In this
image we are interested in terms containing  $duv_{0}$ or $dwuu$.
Let us denote $I_{1}$ (resp. $I_{2}$)  the coefficient at
$duv_{0}$ (resp.  $dwuu$) in $\overline{S}^{12} \overline{S}^{23}
((q^{2}uv-vu)dv).$ A straightforward computation shows that
\begin{equation}
\begin{array}{c}
I_{1}=(q^{3}+q)\alpha _{1}^{\prime \prime }I_{0}+(q+q^{-1})\alpha
_{1}\alpha _{1}^{\prime }[q^{6}-(q+q^{-1})\beta ^{\prime
}+q^{2}\gamma _{1}^{\prime
}]c\delta , \\
I_{2}=(q+q^{-1})I_{0}+(q^{3}+q)\alpha _{1}[\beta ^{\prime
}+q^{6}\gamma
_{1}^{\prime }]c\delta\quad  {\rm where} \\
I_{0}=(q^{3}+q)\alpha _{1}\alpha _{1}^{\prime }\delta ^{2}+\alpha
_{1}[2q^{-3}(1+q^{2})\alpha _{1}^{\prime }+(\beta
_{1}-q^{8}(1+q^{2})\gamma _{1}^{\prime })-1]c\delta .
\end{array}
\end{equation}

Since the element $duv_{0}$\ is that of the highest weight and
the element $(q^{2}uv-vu)dv$  is not,  the coefficient  $I_{1}$
is trivial. Moreover, the condition $\overline{S}(V_{1}\otimes
V^{\prime })\subset (V^{\prime }\otimes V_{1})$ implies that the
coefficient $I_{2}$\ is trivial. These two relations imply $\
c\delta =0.$

Since $c\not=0$ (unless the operator $\overline{S}$ is not
invertible) we have $\delta =0$.

Thus, if we admit the first condition of \r{tri} the only
transposition $\overline{S}=cS$ could preserve the defining ideal
of $\kSq$. But in fact only factors $c=\pm 1$ are compatible with
the centrality of the element $v_0$. This completes the proof.

Now we  pass to showing that even the operator $S$ does not
preserve the flatness of the deformation (for the operators $-S$
and $\pm S^{-1}$ the proof is analogous).

\begin{proposition}
The image of $M_{l}$ w.r.t. the YB operator $S$ \ does not belong
to $M_{r}.$
\end{proposition}

{\bf Proof} By definition of the submodules $M_{l}$ and $M_{r}$,
they consist respectively  of the following elements
\begin{equation}
\mu ^{12}(fdv_{0})\quad {\rm and }\quad \mu^{23}
(d\overline{v}_{0}f), \,\, \forall \,f\in k(S_{q}^{2}).
\end{equation}

Let us show that there exists an element $f\in k(S_{q}^{2})$\
such that $S(\mu ^{12}(fdv_{0}))\notin $\ $M_{r}$. Let $f=u$.  We
have
\begin{eqnarray*}
S(\mu ^{12}(udv_{0})) &=&\mu ^{23}S^{12}S^{23}(udv_{0}) \\
&=&\mu ^{23}(S^{23})^{-1}S^{23}S^{12}S^{23}(udv_{0}) \\
&=&\mu ^{23}(S^{23})^{-1}(d\overline{v}_{0}u).
\end{eqnarray*}

By using
\begin{equation}
(S^{23})^{-1}(d\overline{v}_{0}.u)=\gamma
^{-1}P_{0}^{23}(d\overline{v}_{0}u)+
\beta^{-1}P_{1}^{23}(d\overline{v}_{0}u)+\alpha
^{-1}P_{2}^{23}(d\overline{v}_{0}u),
\end{equation}
we get
\begin{equation}
\mu ^{23}(S^{23})^{-1}(d\overline{v}_{0}u)=\gamma ^{-1}\mu
^{23}(P_{0}^{23}(d \overline{v}_{0}u))+\alpha ^{-1}\mu
^{23}(P_{2}^{23}(d\overline{v}_{0}u)). \label{nonfl}
\end{equation}

We state that there is no element $g\in\kSq$ such that
$\mu^{23}(d\overline{v}_{0}\, g)$ would be equal to the r.h.s. of
\r{nonfl}. Indeed, it could be only an element of the form
$g=\nu\, u,\,\, \nu\in k$. However, since $\alpha\not=\gamma$ and
the both components in \r{nonfl} are non-trivial we conclude that
no appropriate factor $\nu$ exists. This completes the proof.

\begin{remark}

The statement of this proposition can be generalized to other
quantum algebras related to  non-involutary YB operators. The
crucial property of the operator $S$ used in the proof is the
following one: this operator has more than one distinct
eigenvalues and the corresponding components do not vanish in the
algebra $\kSq$.

However, for certain algebras and certain their two-sided modules
the flatness of deformation is valid. Let us consider for
instance,  the quantum cone (it corresponds to the case $c=0$).
The module $\TSq$ defined as above is a flat deformation of its
classical counterpart. This follows from the fact that in the
corresponding quantum coordinate ring defined by $V_0=0,\,\,
V_1=0$  the only component $V_2\subset V^{\ot 2}$ survives. This
prevents us from the effect used in the proof above.

For the same reason,  in quantum geometry dealing with
non-involutary YB operators $S$ it is not convenient to use this
operator (either any other transposition) in order to define a
product $M_1\ot_A M_2$ of two one-sided $A$-modules assuming $A$
to be a quantum algebra looking like that $\kSq$. However, the
product of two modules can be apparently defined  as the quantum
deformation of the product of their classical counterparts. So,
the notation $M_1\ot_A M_2$ must be regarded in this restricted
sense without any transposition of the elements of $A$ and those
of $M_1$ (or $M_2$). \end{remark}

\section{Generic quantum orbits and modules on them}

The main purpose of this section is to generalize the
construction of the cotangent module on quantum sphere to some
other quantum orbits. Also, we define the tangent modules on
these orbits and discuss the problem of equipping the tangent
module with an action on the quantum coordinate ring in question.
All constructions are done in the framework of one-sided modules
over algebras in question. This allows us to hope that these
modules are flatly deformed objects.

First of all, describe quantum orbits in question. Let us begin
with evoking their quasiclassical counterparts (i.e. the
corresponding Poisson structures).

As was shown in \cite{DGS}, on any orbit (of a semisimple
element) $\cO\subset \gggg^*$ where  $\gggg$ is a simple Lie
algebra there exists a family of Poisson-Lie structures (for the
compact form of the Lie algebra this family is labeled by the
elements of $H^2(\cO)$). Morover,  in this family  there exists a
bracket which is compatible with the Kirillov-Kostant-Souriau
one. A quantization of this particular  Poisson bracket can be
realized (at least in the $\gggg=sl(n)$ case) in terms of the
so-called reflection equation  (RE) algebra. The resulting
algebra can be described as an appropriate quotient of  the RE
algebra. Thus, we get an explicit realization of such an algebra
in the spirit of algebraic geometry by means of a   system of
braided algebraic equations. (As for other Poisson-Lie structures
they can be quantized by means of formal series in the spirit of
deformation quantization. Their description in terms of so-called
Hopf-Galois extension is also known, cf. \cite{Sh}, \cite{HM}.)

Let $\gggg=sl(n)$ and $\gq$ be the same as vector space but
equipped with a $\Uqs$-action which is a deformation of the
adjoint action of $\gggg$ onto itself. Let us extend this action
to the space $\ggq$ by means of the coproduct in $\Uqs$ and
decompose it into a direct sum \beq \ggq=I_+ \oplus I_-
\label{decom} \eeq of two $\Uqs$-invariant subspaces $I_+$ and
$I_-$ so that the corresponding algebras
$$\wedge_{\pm}=\wedge_{\pm}(\gq)=T(\gq)/\{I_{\mp}\}$$ would be
flat deformations of the symmetric $\wedge_{+}(\gggg)$ and
skewsymmetric $\wedge_{-}(\gggg)$ algebras respectively. Since
the space $\ggq$ is not multiplicity free (the component
isomorphic to $\gq$ itself comes twice in the decomposition of
$\ggq$ into a direct sum of irreducible $\Uqs$-modules) it is not
evident that decomposition \r{decom} exists.

Nevertheless, it does exist and can be constructed by means of
the RE algebra mentioned above and of a $\Uqs$-covariant pairing.
Let us recall that by the RE algebra one means the algebra
generated by $n^2$ elements $l_i^j,\,\, 1\leq i, \, j\leq n$
subject to the relations \beq SL_1SL_1-L_1SL_1S=0 \label{RE} \eeq
where $L=(l_i^j)$ is the matrix with the entries $l_i^j$ and
$L_1=L\ot \id$.

It is known that this algebra has the center generated by the
elements $C^p_q=\tr_q L^p, p=1,..., n$ where $\tr_q$ is the
quantum analogue  of the usual trace. Then the space $\,\,\span
(l_i^j)\,\,$ is a direct sum of a one-dimensional $\Uqs$-module
generated by $\tr_q L$ and a $n^2-1$-dimensional one which we
identify with $\gq$. Then the space $I_-$ can be treated as the
l.h.s. of \r{RE} modulo the elements of the form \beq l\ot\tr_q
L,\,\, \tr_q L\ot l, \,\,l\in \span (l_i^j). \label{modulo} \eeq
In virtue of \cite{L} the algebra $\wp$ is a flat deformation of
its classical counterpart. (Let us note that the RE algebra and
adjacent objects are also well defined for non-quasiclassical YB
operators, cf. \cite{GPS}).

We introduce the space $I_+$ as that orthogonal to $I_-$  w.r.t.
the pairing $$ (\,\,,\,\,):\ggq\ot \ggq\to k,\,\,
(\,\,,\,\,)=<\,\,,\,\,>\,<\,\,,\,\,>^{23} $$ where
$<\,\,,\,\,>:\ggq\to k$ is a (unique up to a factor)
$\Uqs$-covariant pairing. Then  following \cite{D} we can state
that the algebra  $\wm$ is a flatly deformed object as well.

Now, let us introduce "the generic quantum orbits" by the
following system of equations \beq C^p_q-c^p=0, \,\, c^p\in
k,\,\, p=1,..., n. \label{ideal} \eeq

The constant $c^1$ is equal to 0 while the other constants are
assumed to be generic. Let us note $k(M_q)$ the quotient of the
RE algebra over the ideal $\{J\}$ generated by the l.h.s.
elements of \r{ideal}. This "quantum coordinate ring" is a   flat
deformation of coordinate ring of a generic orbit in $sl(n)^*$.

As for q-analogues of other orbits of semisimple elements in
$sl(n)^*$ the reader is referred to \cite{DGK} where the case of
the "$\bbbc\bbbp^n$ type orbits" was studied.

Now let us introduce quantum analogues of the cotangent module
and its exterior powers on the orbits in question. (In the sequel
all the modules are left.)

Consider the elements $d\, C^p_q$ looking like $d\,v_0$ of the
previous section. This means that the differential $d$ is applied
only to the last factor of the element $C^p_q$. Let us multiply
the elements $d\, C^p_q$ by those of $\kMq$ from the  left and
the elements of $\wm$ (in the sense of the algebra $\wm$) from
the right.

Now, consider the quotient of the left $\kMq$-module
${\kMq}\ot\wm^l$ over its submodule formed by  the elements $$
f_p\,d\, C^p_q\wedge g_p,\,\, f_p\in {\kMq},\, g_p\in \wm^{l-1} $$
(hereafter $\wm^l$ is the degree $l$ homogeneous component of the
algebra $\wm$). Conjecturally, this quotient is a flat
deformation of the space of degree $l$ differential forms. A
proof of this conjecture for the quantum sphere is given in
\cite{AG}. Also suggested in that paper was a $\Uqs$-covariant de
Rham type complex which was a deformation of its classical
counterpart without making use of the Leibniz rule.

However, the orbits in question are not multiplicity free any
more and the problem of constructing a $\Uqs$-covariant complex
which would be a flat deformation of the usual de Rham one
becomes more delicate. Nevertheless, if we are only interested in
q-analogues of 2-forms being generators of the cohomology ring on
quantum orbits in question we can explicitly construct them in
the following way (by analogy with the classical case).

It suffices to apply the q-cobracket to the last factor of each
element $C_q^p$ and treat its image as an element of $\wedge^2_-$
(i.e. by realizing it as a sum of the summands $d\,x_i\wedge
d\,x_j$). By definition, the q-cobracket is  the inverse (in some
natural sense) of the q-Lie bracket whose construction was
featured in \cite{LS}.

Let us call the above quotient module corresponding to the case
$l=1$ as {\em the cotangent module} and denote it as $T^*(M_q)$.
Emphasize once more that this module is introduced as a one-sided
module (namely, the left one but in the same way we can introduce
the right one). Moreover, it is introduced explicitly by a system
of equations. In a    similar way we can realize the other
modules defined above.

Now let us pass to defining the tangent module $T(M_q)$ on the
orbits in question. In the classical case the tangent module
$T(M)$ on a given regular affine variety $M$ has an operator
realization by vector fields, i.e., there exists a   map \beq
T(M)\ot k(M) \to k(M)\label{anch} \eeq which commutes with the
$k(M)$-module structure product $$k(M)\ot T(M) \to T(M).$$
Moreover, if $M$ is an orbit in $\gggg^*$ there exists an
embedding $$\gggg\hookrightarrow T(M)$$ such that map \r{anch}
realizes a representation of Lie algebra  $\gggg$ by vector
fields in the coordinate ring $k(M)$.

As for the tangent modules on quantum orbits $k(M_q)$ we define
them by the same system as the cotangent ones. This is motivated
by the fact that in the classical case the tangent and cotangent
modules on  orbits in question are isomorphic. So, conjecturally
the deformation of the tangent module is flat.

However, in the quantum case there exists the problem of an
operator meaning of the tangent module. For the case of the
quantum sphere this problem was solved in  \cite{A}. Namely, it
was shown that for the tangent quantum module $T(S^2_q)$ there
exists a map \beq T(S^2_q) \ot k(S^2_q)\to k(S^2_q) \label{anchor}
\eeq commuting with the module structure product $$ k(S^2_q)\ot
T(S^2_q) \to  T(S^2_q). $$ Also in \cite{A} an embedding was
constructed of the form $$ sl(2)_q\hookrightarrow T(S^2_q) $$
such that map \r{anchor} realized a representation of the q-Lie
algebra $sl(2)_q$ (this means that the relations between the
generators of $sl(2)_q$ in its enveloping algebra are preserved
under map \r{anchor}).

We call {\em braided vector fields} the elements of $T(S^2_q)$
realized as operators on $k(S^2_q)$ via the map \r{anchor}.
However the problem of a similar treatement of the tangent
modules on the quantum orbits in question is still open.

Let us complete the paper with the following remark. There exists
a lot of articles devoted to different aspects of "braided
geometry". However, often they do not consider the problem of
flatness of quantum deformation. Nevertheless, a flat deformation
is rather subtle phenomenon. Even if the flatness is fulfilled for
a  deformation of complexes related to a vector space, in general
it disappears if one tries to restrict the differential algebras
to a "quantum variety". We are sure that the approach making use
of one-sided modules on quantum varieties developed here (as well
as in \cite{AG}, \cite{GS}) is more adequate for the needs of
"braided geometry" on quantum varieties since conjecturally it
allows us to preserve the flatness of deformation (at the expense
of the Leibniz rule).

\end{document}